\documentclass [12pt,a4paper,reqno]{amsart}
\input amssymb.sty
\usepackage{amsmath,amsthm}
\usepackage{amsfonts}
\usepackage{amssymb}
\allowdisplaybreaks

\newcommand{\be}{\begin{equation}}
\newcommand{\ef}{\end{equation}}
\chardef\bslash=`\\ % p. 424, TeXbook
\newtheorem{thm}{Theorem}[section]
\newtheorem*{thm*}{Theorem}

\newtheorem{lem}[thm]{Lemma}
\newtheorem{prop}[thm]{Proposition}

\theoremstyle{definition}

\newtheorem*{remark*}{Remarks}
\newtheorem*{defn*}{Definition}

\theoremstyle{remark}

\numberwithin{equation}{section}

\newcommand{\G}{\Gamma}
\newcommand{\wt}{\widetilde}
\newcommand{\wh}{\widehat}
\newcommand{\fc}{\frac}
\newcommand{\bk}{\bigskip}
\newcommand{\iy}{\infty}
 \renewcommand{\sectionmark}[1]{}

\newcommand{\Be}{Beltrami}
\newcommand{\Gr} {Grunsky}
\newcommand{\hol} {holomorphic}
\newcommand{\qc} {quasiconformal}
\newcommand{\sh} {subharmonic}
\newcommand{\psh} {plurisubharmonic}

\newcommand{\ve}{\varepsilon}

\newcommand{\Te} {Teichm\"{u}ller}
\newcommand{\Ko} {Kobayashi}
\newcommand{\Ca} {Carath\'{e}odory}
\newcommand{\uTs} {universal Teichm\"{u}ller space}
\newcommand{\const}{\operatorname{const}}

 \usepackage{amsfonts}
\newcommand{\field}[1]{\mathbb{#1}}
\newcommand{\g}{\gamma}

\newcommand{\D}{\Delta}
\newcommand{\om}{\omega}
\newcommand{\z}{\zeta}
\newcommand{\ov}{\overline}
\newcommand{\vp}{\varphi}
\newcommand{\hC}{\wh{\field{C}}}
\newcommand{\C}{\field{C}}

\newcommand{\B}{\mathbf{B}}
\newcommand{\T}{\mathbf{T}}
\newcommand{\Belt}{\mathbf{Belt}}
\newcommand{\Hol}{\operatorname{Hol}}
\newcommand{\Pot}{\operatorname{Pot}}

\renewcommand{\a} {\alpha}

\newcommand{\ld}{\lambda}
\newcommand{\kp}{\kappa}

\begin{document}

\title{A short geometric proof of the Zalcman and Bieberbach conjectures}
\author{Samuel L. Krushkal}

\begin{abstract}
We show that complex geometric features of \Te \ spaces create
explicitly the extremals of generic homogeneous \hol \ functionals
on univalent functions. In particular this gives proofs of the
well-known Zalcman and Bieberbach conjectures and many new distortion
theorems.
\end{abstract}

\date{\today\hskip4mm({shortZB.tex})}

\maketitle

\bigskip

{\small {\textbf {2010 Mathematics Subject Classification:} Primary:
30C50, 30C75, 30F60, 32Q45; Secondary 30C55, 30C62}

\medskip

}

\bigskip

\markboth{Samuel L. Krushkal}{The Zalcman and Bieberbach conjectures}
\pagestyle{headings}

\section{Introduction}

Our aim is to show that complex geometry of the \uTs \ and \Te \
space of the punctured disk describes explicitly the extremals of
generic homogeneous \hol \ functionals on univalent functions. This
yields, in particular, the proof of the famous Zalcman and
Bieberbach conjectures and many new distortion theorems.

\subsection{Classes of functions and general homogeneous \hol \
functionals}
The  \hol \ functionals on the classes of univalent functions
depending on the Taylor coefficients of these functions play an
important role in various geometric and physical applications of
complex analysis, for example, in view of their connection with
string theory and with a \hol \ extension of the Virasoro algebra.
These coefficients reflect the fundamental intrinsic features of
conformal maps. Thus estimating them still remains an important
problem in geometric function theory.

We consider the univalent functions on the unit disk $\D = \{|z| <
1\}$ normalized by
$$
f(z) = z + \sum\limits_{n=2}^\iy a_n z^n.
$$
These functions form the well-known class $S$. Their inversions
$F_f(z) = 1/f(1/z)$ form the collection $\Sigma$ of univalent
nonvanishing functions $(F(z) \ne 0)$ on the complementary disk
$\D^* = \{z \in \hC = \C \cup \{\iy\}: \ |z| > 1\}$ with expansions
 \be\label{1.1}
F_f(z) = 1/f(1/z)= z + b_0 + b_1 z^{-1} + b_2 z^{-2} + \dots
\end{equation}
Easy computations yield that the coefficients $a_n$ and $b_j$ are
related by
 \be\label{1.2}
b_0 + a_2 = 0, \quad b_n + \sum \limits_{j=1}^{n} b_{n-j} a_{j+1} +
a_{n+2} = 0, \quad n = 1, 2, ... \ ,
\end{equation}
which implies successively the representations of $a_n$ by $b_j$.
One gets
 \be\label{1.3}
a_n = (- 1)^{n-1} b_0^{n-1} - (- 1)^{n-1} (n - 2) b_1 b_0^{n-3} +
\text{lower terms with respect to} \ b_0;
\end{equation}
in particular,
$$
\begin{aligned}
a_2 &= - b_0, \ \  a_3 = - b_1 + b_0^2, \ \
a_4 = - b_2 + 2 b_1 b_0 - b_0^3,  \\
a_5 &= - b_3 + 2 b_2 b_0 + b_1^2 - 3 b_1 b_0^2 + b_0^4, \\
a_6 &= - b_4 + 2 b_3 b_0 + 2 b_2 b_1 - 3 b_2 b_0^2 - 3 b_1^2 b_0
+ 4 b_1 b_0^3 - b_0^5, \\
a_7 &= b_0^6 - 5 b_1 b_0^4  - b_1 ^3 + 4 b_2 b_0^3 + b_2^2
+ (6 b_1^2 - 3 b_3) b_0^2 \\
&+ 2 b_1 b_3 + (-6 b_1 b_2 + 2 b_4) b_0 - b_5, \ \dots
\end{aligned}
$$
We shall essentially use this connection.

Consider a general \hol \  distortion functional on $S$ of the form
 \be\label{1.4}
J(f) = J(a_2, \dots \ , a_n; (f^{(\a_1)}(z_1)); \dots \ ;
(f^{(\a_p)}(z_p))),
\end{equation}
where $z_1, \dots \ , z_p$ are the distinct fixed points in $\D
\setminus \{0\}$ with assigned orders $m_1, \ \dots \ , m_p$,
respectively, $(f^{(\a_1)}(z_1)) = f^{\prime\prime}(z_1), \dots,
f^{(m_1)}(z_1); \ (f^{(\a_p)}(z_p)) = f^{\prime\prime}(z_p), \dots,
f^{(m_p)}(z_p)$. Assume that $J$ is a polynomial in all of its
variables.

Substituting the expressions of $a_j$ by $b_m$ from (1.2) and
calculating $f^{(q)}(z_j)$ in terms of $F_f$, one obtains a
polynomial $\wt J(F)$ of the Taylor coefficients $b_0, b_1, \dots \
, b_{n-2}$ and of the corresponding derivatives $F_f^{(q)}(\z_j)$ at
the points $\z_j = 1/z_j \in \D^* \setminus \{\iy\}$, regarded as a
representation of $J(f)$ on the class $\Sigma$. Here $q = 2, \ \dots
\ , m_j, \ j = 1, \ \dots \ , p$.

Assume that the functional (1.4) is homogeneous with a degree $d =
d(J)$ (depending on $n$ and $m_1, \ , \dots, \ , m_p$) with respect
to the homotopy
$$
f(z, t) = t^{-1} f(t z) = z + a_2 t + a_3 t^2 + \dots: \ \D \times
\ov \D \to \C
$$
such that $f(z, 0) \equiv z, \ f(z, 1) = f(z)$ so that
$$
J(f_t) = t^d J(f).
$$
This homotopy is a special case of \hol \ motions with complex
parameter $t$ running over the disk $\D$. The functional $\wt J(F)$
on $\Sigma$ admits a similar homogeneity.

The existence of extremal functions of $J(f)$ and $\wt J(F)$ follows
from compactness of both classes $S$ and $\Sigma$ in the topology of
locally uniform convergence on $\D$ and $\D^*$, respectively.

\subsection{The Bieberbach and Zalcman conjectures}
There were several classical conjectures about the coefficients.
They include the Bieberbach conjecture that in the class $S$  the
coefficients are estimated by $|a_n| \le n$, as well as several
other well-known conjectures that imply the Bieberbach conjecture.
Most of them have been proved by the de Branges theorem \cite {DB}.

In the 1960s, Lawrence Zalcman posed the conjecture that {\em for
any $f \in S$ and all $n \ge 3$,
 \be\label{1.5}
|a_n^2 - a_{2n-1}| \le (n-1)^2,
\end{equation}
with equality only for the Koebe function}
 \be\label{1.6}
\kp_\theta(z) = \fc{z}{(1 - e^{i \theta} z)^2} = z +
\sum\limits_2^\iy n e^{- i(n-1) \theta} z^n, \quad 0 \le \theta \le
2 \pi,
\end{equation}
which maps the unit disk onto the complement of the ray
$$
w = -t e^{-i \theta}, \ \ \fc{1}{4} \le t \le \iy.
$$
This remarkable conjecture also implies the Bieberbach conjecture
and remained an intriguing very difficult open problem for all $n >
6$.

The original aim of Zalcman's conjecture was to prove the Bieberbach
conjecture using the famous Hayman theorem on the asymptotic growth
of coefficients of individual functions, which states that {\em for
each $f \in S$, we have the inequality
$$
\lim\limits_{n \to \iy} \fc{|a_n|}{n} = \a \le 1,
$$
with equality only when $f = \kp_\theta$; here} $\a = \lim\limits_{r
\to 1} (1 - r)^2 \max_{|z|=r} |f(z)|$ (see \cite{Ha}).

Indeed, assuming that $n$ is sufficiently large and estimating
$a_{2n-1}$ in (1.5) by $|a_{2n-1}| \le 2n -1$, one obtains
$$
|a_n|^2 \leq (n -1)^2 + |a_{2n -1}| \leq (n - 1)^2 + 2n - 1 = n^2,
$$
which proves the Bieberbach conjecture for this $n$, and
successively for all preceding coefficients.

It was realized almost immediately that the Zalcman conjecture
implies the Bieberbach conjecture, and in a very simple fashion,
without Hayman's result and without other prior results from the
theory of univalent functions.

Note that the case $n = 2$ is rather simple and somewhat
exceptional. The inequality $|a_2^2 - a_3| \le 1$ is well known, but
in this case there are two extremal functions of different kinds:
the Koebe function $\kp_\theta(z)$ and the odd function
 \be\label{1.7}
\kp_{2,\theta}(z) := \sqrt{\kp_\theta(z^2)} = \sum\limits_{n=0}^\iy
e^{i n \theta} z^{2n+1}.
\end{equation}
The estimate (1.1) was established for $n \le 6$ in \cite{Kr4},
\cite{Kr7} (for $n = 4, 5, 6$ without uniqueness of the extremal
function). In \cite{BT}, \cite{Ma}, this conjecture was proved for
certain special subclasses of $S$.

\section{Main theorems}

\subsection{General theorem}
It is well known that the Koebe function $\kappa_\theta$ is extremal
for many variational problems in the theory of conformal maps
(accordingly, its root transforms
 \be\label{2.1}
\kp_{m,\theta}(z) = \kp_\theta(z^m)^{1/m} = \fc{z}{(1 - e^{i \theta}
z^m)^{2/m}} = z + \fc{2 e^{i \theta}}{m} z^{m+1} + \fc{m - 2}{m^2}
z^{2m+1} + \dots, \quad m = 2, 3, \dots,
\end{equation}
are extremal among the maps with symmetries).

Our first main theorem sheds new light on this phenomenon and
provides a large class of functionals maximized by these functions.

\begin{thm}
Let $J(f)$ be a homogeneous polynomial functional on $S$ of the form
(2.4) whose representation $\wt J(F_f)$ in the class $\Sigma$ does
not contain free terms $c_d b_0^d$ but contains nonzero terms with
the coefficient $b_1$ of inversions $F_f$. Then for all $f \in S$,
we have the sharp bound
 \be\label{2.2}
|J(f)| \le \max_m |J(\kappa_{m,\theta})|,
\end{equation}
and this maximum is attained on some $\kappa_{m_0,\theta} \ (m_0 \ge
1)$. If $J$ has an extremal with
 \be\label{2.3}
b_1 = a_2^2 - a_3 \ne 0,
\end{equation}
then $|b_1| = 1$ and
 \be\label{2.4}
|J(f)| \le \max \{|J(\kappa_\theta)|, |J(\kappa_{2,\theta})|\}.
\end{equation}
\end{thm}

The assumption (2.3) is equivalent to
$$
S_f(0) = - \lim\limits_{z\to \iy} z^4 S_{F_f}(z) \ne 0,
$$
where $S_f$ denotes the {\bf Schwarzian derivative} of $f$ in $\D$
defined by
$$
S_f = (f^{\prime\prime}/f^\prime)^\prime -
(f^{\prime\prime}/f^\prime)^2/2.
$$

The examples of some well-known functionals, for example, $J(f) =
a_2^2 - \a a_3$ with $0 < \a < 1$ and $J(F_f) = b_m \ (m > 1)$, show
that the assumptions on the initial coefficients $b_0$ and $b_1$
cannot be omitted.

\subsection{Applications}
The Zalcman functional
$$
J_n(f) = a_n^2 - a_{2n-1}
$$
is a special case of (2.4) with homogeneity degree $2n - 2$. For
this functional, we obtain from Theorem 2.1 a complete result
proving the Zalcman conjecture.

\bk
\begin{thm} For all $f \in S$ and any $n \ge 3$,
we have the sharp estimate (1.5), with equality only for $f =
\kappa_\theta$.
\end{thm}

As a consequence, one obtains also a new proof of the Bieberbach
conjecture.

Theorem 2.1 also provides other new distortion theorems 
concerning the higher coefficients. These results are presented in
Section 6. In the last section, we show how the proof of Theorem 2.1
yields asymptotic estimating the growth rate of generic homogeneous
functionals on an individual function with \qc \ extension.

\subsection{Connection with geometry of \Te \ spaces}
Our approach to these problems is geometric. Its origins go back to
\cite{Kr7} where the proof of Zalcmann's conjecture for the
initial coefficients was given. 

It suffices to find the bound of $J$ on functions $f$ admitting \qc
\ extensions across the unit circle and close this set in weak
topology determined by locally uniform convergence on $\D$. Denote
the subset of such $f$ by $S^0$ and the set of corresponding $F_f
\in \Sigma$ by $\Sigma^0$.

Such functions are naturally connected with the \uTs \ $\T = \T(\D)$
and the \Te \ space $\T_1 = \T(\D \setminus \{0\})$ of the punctured
disk. Accordingly, the original functional $J(f)$ is lifted to a
\hol \ functional on $\T_1$, and its sharp upper bound is obtained
using deep geometric features of this space. Application of metrics
of negative generalized curvature in the lines of \cite{Kr7} allows
us to estimate the functional from below giving the same asymptotic
bound.

In fact, we establish that every 
homogeneous \hol \ functional on $S$ satisfying the 
assumptions of Theorem 2.1 determines a complex geodesic in  
the space $\T_1$ generated by some $\kappa_{m,\theta}$.

\section{Background}

We briefly present here certain results underlying the proof of the 
key Theorem 2.1. The exposition is adapted to our special cases.

\subsection{A glimpse at complex geometry of \Te \ spaces $\T$ and
$\T_1$}

\bk\noindent{\em (a)} \ First recall that the \uTs \ $\T$ is the
space of quasisymmetric homeomorphisms of the unit circle $S^1 =
\partial \D$ factorized by M\"{o}bius maps. All \Te \ spaces have 
their isometric copies in $\T$. 

The canonical complex
Banach structure on $\T$ is defined by factorization of the ball of
the {\bf \Be \ coefficients} (or complex dilatations)
 \be\label{3.1}
\Belt(\D)_1 = \{\mu \in L_\iy(\C): \ \mu|\D^* = 0, \ \|\mu\| < 1\},
\end{equation}
letting $\mu_1, \mu_2 \in \Belt(\D)_1$ be equivalent if the
corresponding \qc \ maps $w^{\mu_1}, w^{\mu_2}$ (solutions to the
\Be \ equation $\partial_{\ov{z}} w = \mu \partial_z w$ with $\mu =
\mu_1, \mu_2$) coincide on the unit circle $S^1 = \partial \D^*$
(hence, on $\ov{\D^*}$). The equivalence classes $[w^\mu]$ are in
one-to-one correspondence with the Schwarzian derivatives $S_w$ of
$w = F^\mu$ on $\D^*$.

The smallest dilatation $k(w) = \inf \|\mu_w\|_\iy$ among \qc \
extensions of univalent $w|\D^* \in \Sigma^0$ onto $\hC$ is called
the \textbf{\Te \ norm} of $w$.

Note that for each locally univalent function $w(z)$ on a simply
connected hyperbolic domain $D \subset \hC$, its Schwarzian
derivative $S_w$ belongs to the complex Banach space $\B(D)$ of
hyperbolically bounded \hol \ functions on $D$ with the norm
$$
\|\vp\|_\B = \sup_D \ld_D^{-2}(z) |\vp(z)|,
$$
where $\ld_D(z) |dz|$ is the hyperbolic metric on $D$ of Gaussian
curvature $- 4$; hence $\vp(z) = O(z^{-4})$ as $z \to \iy$ if $\iy
\in D$. In particular, for $D = \D$,
 \be\label{3.2}
\ld_\D(z) = 1/(1 - |z|^2).
\end{equation}
The space $\B(D)$ is dual to the Bergman space $A_1(D)$, a subspace
of $L_1(D)$ formed by integrable \hol \ functions on $D$.

The Schwarzians $S_{w^\mu}(z)$ with $\mu \in \Belt(\D)_1$ range over
a bounded domain in the space $\B = \B(\D^*)$. This domain models
the \uTs \ $\T$, and the factorizing projection
$$
\phi_\T(\mu) = S_{w^\mu}: \ \Belt(\D)_1 \to \T
$$
is a \hol \ map from $L_\iy(\D)$ to $\B$. This map is a split
submersion, which means that $\phi_\T$ has local \hol \ sections
(see, e.g., [GL]).

Both equations $S_w = \vp$ and $\partial_{\ov z} w = \mu
\partial_z w$ (on $\D^*$ and $\D$, respectively) determine their
solutions in $\Sigma^0$ up to translations $w \mapsto w + b_0$. To
determine a solution  $w^\mu$ uniquely, we add the condition
$w^\mu(0) = 0$ going over from $w^\mu$ to the maps
$$
w_1^\mu(z) = w^\mu(z) - w^\mu(0) = z - \fc{1}{\pi} \iint\limits_\D
\fc{\partial w^\mu}{\partial \ov \z} \left( \fc{1}{\z - z} -
\fc{1}{\z} \right) d \xi d \eta \quad (\z = \xi + i \eta).
$$
Then the values $w^\mu(z_0)$ (for any fixed $z_0 \in \C$) and the
Taylor coefficients $b_1, b_2, \dots$ of $w^\mu \in \Sigma^0$ depend
holomorphically on $\mu \in \Belt(\D)_1$ and on $S_{w^\mu} \in \T$.

The points of \Te \ space $\T_1 = \T(\D^0)$ of the punctured disk
$\D^0 = \D \setminus \{0\}$ are the equivalence classes of \Be \
coefficients $\mu \in \Belt(\D)_1$ so that the corresponding \qc \
automorphisms $w^\mu$ of the unit disk coincide on both boundary
components (unit circle $S^1 = \{|z| =1\}$ and the puncture $z = 0$)
and are homotopic on $\D \setminus \{0\}$. This space can be endowed
with a canonical complex structure of a complex Banach manifold and
embedded into $\T$ using uniformization.

Namely, the disk $\D^0$ is conformally equivalent to the factor
$\D/\G$, where $\G$ is a cyclic parabolic Fuchsian group acting
discontinuously on $\D$ and $\D^*$. The functions $\mu \in
L_\iy(\D)$ are lifted to $\D$ as the \Be \ $(-1, 1)$-measurable
forms  $\wt \mu d\ov{z}/dz$ in $\D$ with respect to $\G$, i.e., via
$(\wt \mu \circ \g) \ov{\g^\prime}/\g^\prime = \wt \mu, \ \g \in
\G$, forming the Banach space $L_\iy(\D, \G)$.

We extend these $\wt \mu$ by zero to $\D^*$ and consider the unit ball
$\Belt(\D, \G)_1$ of $L_\iy(\D, \G)$. Then the corresponding
Schwarzians $S_{w^{\wt \mu}|\D^*}$ belong to $\T$. Moreover, $\T_1$
is canonically isomorphic to the subspace $\T(\G) = \T \cap \B(\G)$,
where $\B(\G)$ consists of elements $\vp \in \B$ satisfying $(\vp
\circ \g) (\g^\prime)^2 = \vp$ in $\D^*$ for all $\g \in \G$. Most
of the results about the \uTs \ presented in Section 1 extend
straightforwardly to $\T_1$.

Due to the Bers isomorphism theorem, the space $\T_1$ is
biholomorphically equivalent to the Bers fiber space
$$
\mathcal F(\T) = \{\phi_\T(\mu), z) \in \T \times \C: \ \mu \in
\Belt(\D)_1, \ z \in w^\mu(\D)\}
$$
over the \uTs \ with \hol \ projection $\pi(\psi, z) = \psi$ (see
\cite{Be2}). This fiber space is a bounded domain in $\B \times \C$.

We shall denote the equivalence classes of $\mu \in \Belt(\D)_1$ in
$\T$ and $\T_1$ by $[\mu]$ and $[\mu]_1$ (also by $[w^\mu]$ and
$[w^\mu]_1$), respectively.

\bk 
Let $\wt \T$ denote one of the spaces $\T, \ \T_1$. It is a
complex Banach manifold, thus it possesses the invariant \Ca \ and \Ko
\ distances (the smallest and the largest among all holomorphically
non-expanding metrics). Denote these metrics by $c_{\wt \T}$ and
$d_{\wt \T}$, and let $\tau_{\wt \T}$ be the intrinsic \Te \ metric
of this space canonically determined by \qc \ maps. The
corresponding differential (infinitesimal) Finsler forms of these
metrics are defined on the tangent bundle $\mathcal T \wt \T$ of
$\wt \T$. Then $c_{\wt \T}(\cdot, \cdot) \le d_{\wt \T}(\cdot,
\cdot) \le \tau_{\wt \T}(\cdot, \cdot)$, and by the Royden-Gardiner
theorem the \Ko \ and \Te \ metrics (and their infinitesimal forms)
are equal, see, e.g. \cite{EKK}, \cite{EM}, \cite{GL}, \cite{Ro1}.

\bk\noindent{\em (b)} \ For the spaces $\T$ and $\T_1$, there is a
much stronger result established in \cite{Kr6}, \cite{Kr9}.

\bk\noindent{\bf Theorem A}. {\em The \Ca \ metric of the space 
$\wt \T$ coincides with
its \Ko \ metric, hence all invariant non-expanding metrics on $\wt
\T$ are equal its \Te \ metric, and
 \be\label{3.3}
 c_{\wt \T}(\vp, \psi) = d_{\wt \T}(\vp,\psi) =
 \tau_{\wt \T}(\vp, \psi) = \inf \{
d_\D(h^{-1}(\vp), h^{-1}(\psi)): \ h \in \Hol(\D, \wt \T)\},
\end{equation}
where $d_\D$ denotes the hyperbolic metric of the unit disk of
curvature $-4$ (i.e., with the differential form (3.2)).

Similarly, the infinitesimal forms of these metrics coincide with
the Finsler form of $\tau_{\wt \T}$ and have \hol \ sectional
curvature $-4$.}

\bk
Such a theorem has been proved in \cite{Kr6} for the \uTs. This
proof is complicated and involves the technique of the \Gr \
coefficient inequalities. A much simpler proof was given recently in
\cite{Kr9}, and the same arguments work for $\wt \T = \T_1$.

Theorem A is one of the main ingredients in the proof of our main
theorems. It also has many other applications.
In particular, the \Te \ extremal disks are simultaneously geodesic
for all non-expanding invariant metrics (cf. \cite{EKK}, \cite{Ve}).

Combining Theorem A with Golusin's improvement of Schwarz's lemma, one
derives the following sharp estimate of the growth of \hol \ maps on
geodesic disks.

\begin{prop} \cite{Kr9}  If the restriction of a \hol \ map 
$h: \ \wt \T \to \D$ onto a geodesic disk
$\D(\mu_0) = \{\phi_{\wt \T}(t \mu_0/\|\mu_0\|_\iy): \ |t| < 1\}$  
has at the origin zero of order $m \ge 1$, i.e.,
$$
h_{\mu_0}(t) := h \circ \phi_{\wt \T}(t \mu_0/\|\mu_0\|_\iy) = c_m t^m 
+ c_{m+1} t^{m+1} + \dots,
$$
then the growth of $|h|$ on this disk is estimated by
 \be\label{3.4}
\begin{aligned} 
|h_{\mu_0}(t)| &\le |t|^m (|t| + |c_m|)/(1 + |c_m||t|)  \\
&= \tanh  d_{\wt \T} \Bigl(\mathbf 0, \phi_{\wt \T} 
\Bigl(|t|^m \ \fc{|t| + |c_m|}{1 + |c_m||t|}\fc{\mu_0}{\|\mu_0\|_\iy}
\Bigr)\Bigr) \le \tanh  d_{\wt \T} \Bigl(\mathbf 0, \phi_{\wt \T} 
\Bigl(t^m \fc{\mu_0}{\|\mu_0\|_\iy}\Bigr)\Bigr). 
\end{aligned}
\end{equation}
The equality in the right inequality occurs (even for one $t_0 \ne
0$) only when $|c_p| = 1$; then $h_{\mu_0}(t)$ is a hyperbolic
isometry of the unit disk and all terms in (3.4) are equal.
\end{prop}

Golusin's version mentioned above asserts that a \hol \ function
$$
g(t) = c_m t^m + c_{m+1} t^{m+1} + \dots: \D \to \D \quad (c_m \neq
0, \ \ m \ge 1)
$$
is estimated in $\D$ by
$$
|g(t)| \le |t|^m \fc{|t| + |c_m|}{1 + |c_m| |t|},
$$
and the equality occurs only for $g_0(t) = t^m (t+ c_m)/(1 + \ov c_m
t)$ (see [Go, Ch. 8]). 

On the other hand, it follows from Theorem A and weak$^*$ compactness 
of the closure of $\wt \T$ in $\B \times \D$ that for any fixed 
$t_0 \ne 0$ there is a \hol \ map $j(\vp): \ \wt \T \to \D$ such that
$$
d_\D(0, j \circ \phi_{\wt \T}(t_0 \mu_0^*)) = 
c_{\wt \T} (\mathbf 0, \phi_{\wt \T} (t_0  \mu_0^*)) 
= d_{\wt \T} (\mathbf 0, \phi_{\wt \T} (t_0 \mu_0^*)).
$$ 
where $\mu_0^* =  \mu_0/\|\mu_0\|_\iy$. 
Thus, letting 
$$
\eta(t) = |t|^m (|t|+ |c_m|)/(1 + |c_m| |t|) \le |t|, 
$$
one derives
$$
|h_{\mu_0}(t_0)| \le j \circ \phi_{\wt \T}(\eta(t_0) \mu_0^*) 
= \tanh d_{\wt \T}(\mathbf 0,
 \phi_{\wt \T}(\eta(t_0) \mu_0^*)) 
\le \tanh  d_{\wt \T}(\mathbf 0, \phi_{\wt \T}(|t_0| \mu_0^*))
$$
which implies (3.4) (for details see \cite{Kr9}). 

There is also a differential analog of the inequalities (3.4) which
involves the infinitesimal metrics $\mathcal C_{\wt \T}$ and $\mathcal
K_{\wt \T}$. It will not be used here.

\subsection{A holomorphic homotopy of univalent function}
Similar to the functions in $S$, we define for each $F \in \Sigma$
with expansion (1.1) the complex homotopy
 \be\label{3.5}
F_t (z) = t F \left( \fc{z}{t} \right) =  z + b_0 t + b_1 t^2 z^{-1}
+ b_2 t^3 z^{-2} + ...: \ \D^* \times \D \to \hC
\end{equation}
so that $F_0(z) \equiv z$. Then $S_{F_t}(z) = t^{-2} S_F(t^{-1} z)$,
and moreover, this point-wise map determines a \hol \ map
 \be\label{3.6}
h_F(t) =  S_{F_t}(\cdot): \ \D \to \B
\end{equation}
(see, e.g. \cite{Kr3}). This map generates the {\bf homotopy disks}
$\D(S_F) = h_F(\D)$ of $F$ in the space $\T$  
and its covers in $\T_1$ defined by 
$$ 
\D_1(S_F, b_0) := \{(S_{F_t}, t b_0): \ |t| < 1\}.  
$$  
These disks are \hol \ at noncritical points of map (3.6) and 
foliate both spaces and the set $\Sigma^0$. 

The dilatations of the homotopy maps are estimated by

\begin{prop} \cite{Kr3} (a) \ Each homotopy map $F_t$
of $F \in \Sigma$ admits $k$-\qc \ extension to the complex sphere
$\hC$ with $k \le |t|^2$. The bound $k(F_t) \le |t|^2$ is sharp and
occurs only for the maps
  \be\label{3.7}
F_{b_0,b_1} (z) = z + b_0 + b_1 z^{-1}, \quad |b_1| = 1,
\end{equation}
whose homotopy maps
 \be\label{3.8}
F_{b_0,b_1t^2}(z) = z + b_0 t + b_1 t^2 z^{-1}
\end{equation}
have the affine extensions $\wh F_{b_0,b_1t^2}(z) = z + b_0 t + b_1
t^2 \ov z$ onto $\D$.

(b) \ If $F(z) = z + b_0 + b_m z^{-m} + b_{m+1} z^{-(m+1)} + \dots \
(b_m \ne 0)$ for some integer $m > 1$, then the minimal dilatation
of extensions of $F_t$ is estimated by $k(F_t) \le |t|^{m+1}$; this
bound also is sharp.
\end{prop}

In the second case,
$$
h_F(0) = h_F^\prime(0) = \dots = h_F^{(m)}(0) = \mathbf 0, \
h_F^{(m+1)}(0) \ne \mathbf 0,
$$
and due to \cite{KK},
 \be\label{3.9}
k(F_t) = \fc{m + 1}{2} |b_m| |t|^{m+1} + O(t^{m+2}), \quad t \to 0.
\end{equation}
This bound is sharp; it holds for the maps
 \be\label{3.10}
 F_{m,t}(z) = \fc{1}{\kp_{m,t}(1/z)} = z \left(1 -
\fc{t}{z^{m+1}}\right)^{2/(m+1)} = z - \fc{2 t}{m + 1} \fc{1}{z^m} +
\dots,  \quad |t| \le 1,
\end{equation}
whose extremal extension to $\C$ has \Be \ coefficient
$\mu_{F_{m,t}}(z) = t |z|^{m-1}/z^{m-1}$ for $|z| < 1$.

\bk The \Te \ geodesic (extremal) disks
$$
\D(\psi) = \{\phi_{\wt \T}(t \mu_0): \ t \in \D\}
$$
in the spaces $\T$ and $\T_1$ are generated by $F \in \Sigma^0$
having extremal extensions with \Be \ coefficients $\mu_t(z) = t
|\psi(z)|/\psi(z)$, where $\psi$ is a \hol \ integrable quadratic
differential on $\D$ and $\D \setminus \{0\}$, respectively. Such 
an extension is unique (up to a constant factor of $\psi$) in their
equivalence classes.

In particular, any homotopy function $F_t$ has such an extremal
extension. The \Te \ disks foliate dense subsets in $\T$ and $\T_1$
(and in $\Sigma^0$); cf. \cite{GL}, \cite{St}.

\subsection{Two generalizations of Gaussian curvature and circularly
symmetric metrics}

The proof of Theorem 2.1 involves also \sh \ conformal metrics
$\ld(t) |dt|$ on the disk (with $\ld(t) \ge 0$) having the curvature
at most $- 4$ in a somewhat generalized sense. As is well-known, the
Gaussian curvature of a $C^2$-smooth metric $\ld > 0$ is defined by
$$
\kappa_\ld = - \fc{\D \log \ld}{\ld^2},
$$
where $\D$ means the Laplacian $4 \partial \ov{\partial}$.

A metric $\ld(t) |dt|$ in a domain $D$ on $\C$ (or on a Riemann
surface) has curvature less than or equal to $K$ \textbf{in the
supporting sense} if for each $K^\prime > K$ and each $t_0$ with
$\ld(t_0) > 0$, there is a $C^2$-smooth \textbf{supporting metric}
$\wh \ld$ for $\ld$ at $t_0$ (i.e., such that $\wh \ld(t_0) =
\ld(t_0)$ and $\wh \ld(t) \le \ld(t)$ in a neighborhood of $t_0$)
with $\kappa_{\wh \ld}(t_0) \le K^\prime$, or equivalently,
 \be\label{3.11}
\D \log \ld \ge - K \ld^2.
\end{equation}
A metric $\ld$ has curvature at most $K$ \textbf{in the potential
sense} at $z_0$ if there is a disk $U$ about $t_0$ in which the
function
$$
\log \ld + K \Pot_U(\ld^2),
$$
where $\Pot_U$ denotes the logarithmic potential
$$
\Pot_U h = \fc{1}{2 \pi} \int\limits_U h(\z) \log |\z - t| d \xi d
\eta \quad (\z = \xi + i \eta),
$$
is \sh. Since $\D \Pot_U h = h$ (in the sense of distributions), one
can replace $U$ by any open subset $V \subset U$, because the
function $\Pot_U(\ld^2) - \Pot_V(\ld^2)$ is harmonic on $U$. The
inequality (3.11) holds for the generic \sh \ metrics also in the
sense of distributions. Note also that the condition of having
curvature at most $- K$ in the potential sense is invariant under
conformal maps.

Due to \cite{Ro2}, a conformal metric of curvature at most $K$ in
the supporting sense has curvature at most $K$ also in the potential
sense.

The following lemma concerns  the circularly symmetric (radial)
metrics and is a slight improvement of the corresponding Royden's
lemma \cite{Ro2} to  singular metrics with a prescribed singularity
at the origin.

\begin{lem} \cite{Kr7} Let $\ld(|t|) d |t|$ be a circularly symmetric
\sh \ metric on $\D$ such that
 \be\label{3.12}
\ld(r) = m c r^{m-1} + O(r^{m}) \quad \text{as} \ \ r \to 0 \ \
\text{with} \ \ 0 < c \le 1 \ \ (m = 1, 2, \dots), 
\end{equation}
and this metric has curvature at most $- 4$ in the potential sense.
Then
$$
\ld(r) \ge \fc{m c r^{m-1}}{1 - c^2 r^{2m}}.
$$
\end{lem}

Note that all metrics subject to (3.12) are dominated by $\ld_m(t) =
m |t|^{m-1}/(1 - |t|^{2m})$.

\bk

\section{Proof of Theorem 2.1} 

\noindent
$\mathbf 1^0$. \ One may assume that the degree $d$ of $J$ is even, 
replacing, if needed, this functional by its square $J^2$.

Using the relations (1.2), we represent $J$ as a polynomial
functional on $\Sigma$, which takes the form
 \be\label{4.1}
J(f) = \wt J(F_f) = \wt J(b_0, b_1, \dots, b_{2n-3};
F_f^{\prime\prime}(\z_1), \dots, F_f^{(m_1)}(\z_1); \dots,
F_f^{\prime\prime}(\z_p), \dots, F_f^{(m_p)}(\z_p)),
\end{equation}
where $b_0 = - a_2$ and $\z_j = 1/z_j$. As was mentioned above, the
admissible values of $b_0$ for $F(z) = z + b_0 + b_1 z^{-1} + \dots
\in \Sigma^0$ with $F(0) = 0$ range over the closed domain
$F(\ov{\D}) = \hC \setminus F(\D^*)$.

The functional $\wt J(F) = J(f)$ on $F = F_f \in \Sigma^0$ extends
to a \hol \ functional $\mathcal J$ on the fiber space $\mathcal
F(\T) = \{(S_F, b_0)\}$ and thereby on the space $\T_1$, letting 
for $\mu \in
\Belt(\D)_1$ and $f^{\wt \mu} \in S^0$ with $\wt \mu(z) = \mu(1/z)
z^2/\ov z^2$,
 \be\label{4.2}
\mathcal J(S_{F^\mu}, b_0(F^\mu)) = J(f^{\wt \mu}) \quad
(b_0(F^\mu)) = - a_2(f^{\wt \mu})).
\end{equation}
We rescale $\mathcal J$ by
$$
\mathcal J^0(S_F, b_0) = \fc{\mathcal J(S_F, b_0)}{M(J)} \quad
\text{with} \ \ M(J) = \max_S |J(f)|
$$
to have a \hol \ map of $\mathcal F(\T)$ to the unit disk. For any
fixed point $z_{*} \in \ov D$, one determines by (4.1) a \hol \
function
$$
g_{*}(\vp) = \mathcal J^0(- F(z_{*}), \{b(\vp)\},
\{F^{(m_j)}(\z_j(\vp))\}) : \ \T \to \D
$$
where $\{b(\vp)\}$ and $\{F^{(m_j)}(\z_j(\vp))\}$ denote the
collections $(b_1, \dots \ , b_{2n-3})$ and
$$
(F^{\prime\prime}(\z_1), \dots, F^{(m_1)}(\z_1); \dots; \
F^{\prime\prime}(\z_p), \dots, F^{(m_p)}(\z_p)),
$$
respectively, regarded as functions of the Schwarzians $\vp = S_F
\in \T$.

In view of the maximum principle, it suffices to use only the
boundary points $z_{*}$. We select on the unit circle $S^1$ a dense
subset
$$
e = \{z_1, z_2, \dots , \ z_m, \dots\}, 
$$
so that the corresponding sequence of \hol \ maps 
 \be\label{4.3}
g_m(\vp) = \wt J^0(- F(z_m), \{b(\vp)\};
\{F^{(m_j)}(\z_j(\vp))\}): \ \T \to \D, \quad m = 1, 2, \dots
\end{equation}
satisfies
  \be\label{4.4}
\sup_m |g_m(S_F)| = \sup_{\T_1} |\mathcal J^0(S_F, a_2)| =
\sup_{\Sigma^0} |\wt J^0(F)| = \max_S |J(f)|/M(J).
\end{equation}

\bk\noindent
$\mathbf 2^0$. \ First suppose that there exists an extremal 
of $J(f)$ satisfying the assumption (2.5), and consider first the
functions $f \in S$ obeying this inequality. The set of the
corresponding Schwarzians $S_{F_f}$ is dense in $\T$, and their maps
(3.6) satisfy
$$
h_F(0) = h_F^\prime(0) = \mathbf 0, \ \ h_F^{\prime\prime}(0) \ne
\mathbf 0.
$$

We split every homotopy function $F_t$ of $F = F_f$ by
$$
F_t(z) = z + b_0 t + b_1 t^2 z^{-1} + b_2 t^3 z^{-2} + \dots =
F_{b_0,b_1t^2}(z) + h(z,t).
$$
For sufficiently small $|t|$, the remainder $h$ is estimated by
$h(z,t) = O(t^3)$ uniformly in $z$ for all $|z| \ge 1$. Then, by the
well-known properties of Schwarzians, we have
$$
S_{F_t}(z) = S_{F_{b_0,b_1t^2}}(z) + \om(z, t), 
$$
where the remainder $\om$ is uniquely determined by the chain rule
$$
S_{w_1\circ w}(z) = (S_{w_1} \circ w) (w^\prime)^2(z) + S_w(z),
$$
and is estimated in the norm of $\B$ by
 \be\label{4.5}
\|\om(\cdot,t)\|_\B = O(t^3), \quad t \to 0;
\end{equation}
 this estimate is uniform for $|t| < t_0$ (cf., e.g. \cite{Be1},
\cite{Kr1}). Hence, in view of holomorphy, every map (4.3) satisfies
for small $|t|$,
$$
g_m(S_{F_t}) = g_m(S_{F_{0,b_1t^2}}) + O(t^{d+1}),
$$
where the term $O(t^{d+1})$ is estimated uniformly for all $n$, and
therefore,
  \be\label{4.6}
\mathcal J(S_{F_t}, b_0 t) = \mathcal J(S_{F_{b_0,b_1t^2}}, b_0) +
O(t^{d+1}), \quad t \to 0.
\end{equation}
Since
$\mathcal J(S_{F_t}, b_0 t) = t^d \mathcal J(S_F, b_0)$ 
(in view of $d$-homogeneity of the functionals $\wt J$ and
$\mathcal J$), we also have 
 \be\label{4.7}
t^d \mathcal J(S_F, b_0) = \mathcal J(S_{F_{b_0,b_1t^2}}, b_0 t) +
O(t^{d+1}),  \quad t \to 0. 
\end{equation}

The values $\mathcal J(S_{F_{0,b_1t^2}}, b_0 t)$ can be sharply
estimated from above by Proposition 3.1. To this end, denote by 
$s$ the canonical complex parameter on the \Te \ disks in 
$\mathcal F(\T)$ generated by admissible (that is, nonvanishing 
on $\D^*$) functions
$$
F_{b_0,s}(z) = z + b_0 + s z^{-1},
$$
whose extremal extensions onto $\ov \D$ are the affine maps
$$
z \mapsto z + b_0 + s \ov z.
$$
All these disks cover the underlying \Te \ disk  $\D(S_{F_0,s})$ in
the base space $\T$ (note that $S_{F_0,s} = S_{F_{b_0,s}}$;  
the functions $F_{0,s}$ with $b_0 = 0$ are associated with  
points of $\T$ in view of normalization $F_f(0) = 0$). 

If $F_{b_0,s}$ is admissible only for $|s| < s_0 < 1$, one can
reparametrize it using the parameter $\sigma = s/s_0$ which runs
over the unit disk. Then, for each $b_0$, the map
$$
\sigma \mapsto (S_{F_{b_0,\sigma}}, b_0 \sigma), \quad \sigma \in
\D,
$$
is a complex geodesics in the space $\mathcal F(\T)$ and
$$
d_{\T_1}(\mathbf 0, (S_{F_{b_0,s}}, b_0 s)) = d_\T(\mathbf 0,
S_{F_{b_0,s}}),
$$
and similarly for the \Ca \ distances.

By (3.10), the parameters $s$ and $t$ are related near the origin by
  \be\label{4.8}
s = b_1 t^2  + O(t^3) \quad (b_1 \ne 0),
\end{equation}
so the restrictions of $\mathcal J^0$ and of $g_m$ to the
indicated \Te \ disks both have at the origin zero of order $d/2$.  
We have 
$$
\mathcal J^0(S_{F_{b_0,s}}) = \beta_{d/2}(b_0) s^{d/2} 
+ \beta_{d/2+1}(b_0) s^{d/2+1} + \dots,    
$$
and after estimating this map by Proposition 3.1 (applied to 
$\wt \T = \T_1$ and $m = d/2$), 
$$
|\mathcal J^0(S_{F_{b_0,s}}| 
\le |s|^{d/2} \fc{|s| + |\beta_{d/2}(b_0)|}{1 + |\beta_{d/2}(b_0)||s|},   
\quad |s| < 1.   
$$
Replacing $s$ by (4.8) and applying the relation (4.6), one obtains 
$$
|\mathcal J^0(S_{F_t}, b_0 t)| = 
|\mathcal J^0(S_{F_{b_0,b_1t^2}}, b_0 t)| + O(t^{d+1}) 
\le |\beta_{d/2}(b_0)| |b_1|^{d/2} |t|^d +  O(t^{d+1}),  
$$
and after maximizing over admissible $b_0$,  
 \be\label{4.9}  
\max_{b_0} |\mathcal J^0(S_{F_t}, b_0 t)| = 
\max_{b_0} |\mathcal J^0(S_{F_{b_0,b_1t^2}}, b_0 t)| + O(t^{d+1})
\le \max_{b_0} |\beta_{d/2}(b_0)| |b_1|^{d/2} |t|^d + O(t^{d+1}), 
 \end{equation} 
where all ratios $O(t^{d+1})/t^{d+1}$ remain uniformly bounded 
as $t \to 0$.

\bk Now our goal is to show that the right-hand side of (4.9) 
yields simultaneously the lower asymptotic bound for  
 $|\mathcal J^0(S_{F_t}, b_0 t)|$ for small $|t|$ (cf. \cite{Kr7}),  
and therefore, the last inequality in (4.10) is reduced to an equality.

\begin{lem} For any $F(z) = z + b_0 + b_1 z^{-1} + \dots \in \Sigma^0$,
we have
 \be\label{4.10}
\max_{b_0} |\mathcal J^0(S_{F_t}, b_0 t)| 
= \max_{b_0} |\mathcal J^0(F_{b_0,b_1t^2})| + O(t^{d+1}) 
\ge \max_{b_0} |\beta_{d/2}(b_0)| |b_1|^{d/2} |t|^d + O(t^{d+1}).
\end{equation}
again taking the maximum over admissible $b_0$.
\end{lem}

\medskip\noindent 
\textbf{Proof.} The homotopy disk of any function
$F_{b_0,b_1}$ is extremal and admits the rotational symmetry. 
Accordingly,  $|\mathcal J(F_{b_0,b_1t^2})|$ is circularly symmetric 
on the disk $\D$. 
Define by (4.3) the corresponding functions  
$$
g_{m,b_1}(t) = \mathcal J^0(S_{F_{b_0,b_1t^2}}, - F_{b_0,b_1t^2}(z_m))  
$$ 
and conformal metrics
$$
\ld_{g_{m,b_1}}(t) = g_{m,b_1}^* \ld_\D(t) 
= \fc{|g_{m,b_1}^\prime(t)|}{1 - |g_{m,b_1}(t)|^2}  
$$
(whose Gaussian curvature equals $- 4$ at noncritical points) 
and take the upper envelopes
$$
\mathcal J^0(t) := \sup_m |g_{m,b_1}(t)|, \quad
\ld_{\mathcal J^0}(t) := \sup_m \ld_{g_{m,b_1}}(t). 
$$
Both envelopes are circularly symmetric continuous and \sh \ on $\D$.  
Note also that (cf. (4.4)), 
\be\label{4.11}
\mathcal J^0(t) \le |\mathcal J^0(S_{F_t}, b_0 t)| 
= \sup_m |g_{m,b_1}(S_{F_{0,b_1t^2}})| 
= |\beta_{d/2}(b_0)| |b_1|^{d/2} |t|^d + O(t^{d+1}).
 \end{equation} 

Since
$$
\tanh^{-1} |g_{m,b_1}(r)| = \int\limits_0^{|g_{m,b_1}(r)|} \fc{|d
t|}{1 - |t|^2} = \int\limits_0^r \ld_{g_{m,b_1}}(t) |d t|,
$$
we have 
 \be\label{4.12}
 \tanh^{-1} |\mathcal J^0(r)| 
= \sup_m \int\limits_0^r \ld_{g_{m,b_1}}(t) |d t|  
= \int\limits_0^r \sup_m
\ld_{g_{m,b_1}}(t) |dt| = \int\limits_0^r \ld_{\mathcal J^0}(t) dt.
\end{equation}
The second equality in (4.12) is obtained by taking a monotone
increasing subsequence of metrics
$$
\ld_1 = \ld_{g_{1,b_1}}, \ \ld_2 
= \max (\ld_{g_{1,b_1}}, \ld_{g_{2,b_1}}), \ 
\ld_3 = \max (\ld_{g_{1,b_1}}, \ld_{g_{2,b_1}}, \ld_{g_{3,b_1}}), 
\ \dots
$$
so that
$$
\lim\limits_{p \to \iy} \ld_p(t) =  \sup_m \ld_{g_{m,b_1}}(t).
$$
Combining (4.12) with (4.4) and (4.7), one gets for small $r$,  
 \be\label{4.13}
\tanh^{-1} [\mathcal J^0(S_{F_1}, b_0r)] 
= \tanh^{-1} [\mathcal J^0(r)] + O(r^{d+1}) 
= \int\limits_0^r \ld_{\mathcal J^0}(t) dt + O(r^{d+1}). 
\end{equation}

Now observe that the metric $\ld_{\mathcal J^0}$ has in a 
neighborhood of any 
$t_0 \in \D$ a supporting metric of curvature $-4$, and therefore   
its curvature in $\D$ in the potential sense is at most $-4$.  
Thus this metric can be estimated from below by
Lemma 3.3 (with $m = d$) which implies the lower bound 
 \be\label{4.14}
\ld_{\mathcal J^0}(r) \ge \fc{d C d r^{d-1}}{1 - C^2 r^{2d}}, 
\quad r < 1, 
\end{equation}
where
$$
C = \max_{b_0} |\beta_{d/2}(b_0)| |b_1|^{d/2}.  
$$
Integrating (4.14) over a small radial segment $[0, r]$, 
one obtains
$$
\int\limits_0^r \ld_{\mathcal J^0}(t) dt \ge \tanh^{-1}
(C r^d) + O(r^{d+1}, \quad r \to 0, 
$$ 
which provides after substitution into (4.13) the desired 
estimate (4.11). 

\bk
Comparison of the relations (4.7) and (4.9)-(4.11) yields 
$$ 
\begin{aligned} 
r^d |\mathcal J^0(S_F, b_0)| &= |\mathcal J^0(S_{F_r}, b_0 r)| 
= \max_{b_0} |\wt J^0(F_{b_0,b_1r^2})| + O(r^{d+1}) \\
&= \max_{b_0} |\beta_{d/2}(b_0)| |b_1|^{d/2} r^d 
+ O(r^{d+1}), \quad r \to 0, 
\end{aligned} 
$$
and letting $r \to 0$,
 \be\label{4.15}
\begin{aligned} 
|\mathcal J^0(S_F, b_0)| &= \max_{b_0} |\wt J^0(F_{b_0,b_1})| =
\max_{b_0} \big\vert\beta_{d/2}(b_0)\big\vert |b_1|^{d/2}  \\ 
&= \max_{b_0} \big\vert\beta_{d/2}(b_0)\big\vert \ 
\left(\fc{|S_f(0)|}{6}\right)^{d/2} 
\le 1. 
\end{aligned} 
\end{equation}

This estimate is established for all $f \in S$. 
Since for any extremal function $f_0(z) = z +
\sum\limits_2^\iy a_n^0 z^n$ of the functional $J$ and its 
inversion  
$F_0(z) = F_{f_0}(z) = z + b_0^0 + b_1^0 z^{-1} + \dots$ 
obeying (2.3) must be  
$$
\fc{|\wt J(S_{F_0}, - a_2^0)|}{M(J)} 
= |\mathcal J^0(S_{F_0}, - a_2^0)| = 1 
$$ 
and $|\beta_{d/2}(b_0)| \le 1$, 
it follows from (4.15) that necessarily  
$\max_{b_0} |\beta_{d/2}(b_0)| = 1$ and 
$$
|b_1^0| = \fc{1}{6} |S_{f_0}(0)| = |(a_2^0)^2 - a_3^0| = 1.
$$
As was mentioned, such equalities can only occur when $f_0$ either
is the Koebe function $\kp_\theta$ or it coincides with the odd
function $\kp_{2,\theta}$ defined by (1.7). In addition, the
extremality of $f_0$ implies
$$
|J(f_0)| = M(J) = \max \{|J(\kp_\theta)|, |J(\kp_{2,\theta})|\}.
$$ 

\bk\noindent
$\mathbf 3^0$. \ The functions $f \in S$ with $S_f(0) = 0$ omitted 
above can
be approximated (in $\B$-norm) by $f$ with $S_f(0) \ne 0$ by
applying special \qc \ deformations of the plane given by the
following lemma from [Kr1, Ch. 4]. This lemma softens the strongest
rigidity of conformal maps.

\begin{lem} In a finitely connected domain $D \subset \hC$, let
there be selected a set $E$ of positive two-dimensional Lebesgue
measure and the distinct finite points $z_1, \dots, z_n$ with
assigned nonnegative integers $\a_1, \dots, \a_n$, respectively, so
that $\a_j = 0$ for $z_j \in E$. Then, for sufficiently small $\ve >
0$ and $\ve \in (0, \ve_0)$, for any given system of numbers
$\{w_{sj}\}, \ s = 0, 1, \dots, \a_j, \ \ j = 1, \dots, n$, such
that $w_{0j} \in D$,
$$
|w_{0j} - z_j| \le \ve, \ \ |w_{1j} - 1| \le \ve, \ \ |w_{sj}| \le
\ve \ \ (s = 2, \dots, \a_j, \ j = 1, \dots, n),
$$
there exists a \qc \ automorphism $h_\ve$ of the domain $D$, which
is conformal on the set $D \setminus E$ and satisfies
$h_\ve^{(s)}(z_j) = w_{sj}$ for all $s = 0, 1, \dots, \a_j$ and $j =
1, \dots, n$, with dilatation $\|\mu_{h_\ve}\|_\iy \le M\ve$. The
constants $\ve_0$ and $M$ depend only on $D, \ E$ and the vectors
$(z_1, \dots, z_n), \ (\a_1, \dots, \a_n)$.

If the boundary $\G$ of domain $D$ is Jordan or belongs to the class
$C^{l,\a}$, where $0 < \a < 1$ and $l \ge 1$, one can take $z_j \in
\G$ with $\a_j = 0$ or $\a_j \le l$, respectively.
\end{lem}

Now, let $f \in S^0$ have coefficients $a_2$ and $a_3$ related by
$a_3 = a_2^2$, i.e., $b_1(f) := b_1(F_f) = 0$. Since $f(\D^*)$ is a
domain, one can take there a set $E$ of positive measure and
construct by Lemma 4.2 for a sequence $\ve_n \to 0$ such variations
$h_n = h_{\ve_n}$ of $f$ that for each $n$,
$$
b_1(h_n \circ f)= b_1(f) + O(\ve_n) \ne 0, \quad |J(h_n \circ f)| =
|J(f)| + O(\ve_n) > |J(f)|.
$$
Since, by the previous step,
$$
|J(h_n \circ f)| \le \max \{|J(\kp_\theta)|, |J(\kp_{2,\theta})|\},
$$
the same estimate will hold also for $f$.

\bk\noindent
$\mathbf 4^0$. \ Finally, consider the case when $J$ has no 
extremals
$f_0$ satisfying (2.3), and hence any extremal inversion $F_{f_0}$
is of the form
 \be\label{4.16}
F(z) = z + b_0 + b_m z^{-m} + b_{m+1} z^{-(m+1)} + \dots \ \quad
(b_m \ne 0; \ |z| > 1)
\end{equation}
with $m > 1$. If $m+1$ does not divide $d = d(J)$, we consider  the
functional $J^{m+1}$, which is $d(m + 1)$-homogeneous; otherwise one
can use $J$.

One can apply to $J^{m+1}$ the above arguments, replacing
$F_{0,b_1t^2}$ by the corresponding function $F_{m,t}(z) + b_0$, 
where $F_{m,t}$ is given by (3.10) and $b_0$ is the same as in (4.16). 
Its Schwarzian relates to $S_{F_t}$ by
$S_{F_t} = S_{F_{m,t}} + O(t^{m+1})$ as $t \to 0$.  
Now   
$$
[J(S_{F_{m,t}}, b_0 t)/M(J)]^{m+1} = \beta_d(b_0) \ t^{d(m+1)} + \dots,   
$$
and after applying the asymptotic estimate (3.9), one obtains instead 
of (4.15) the bound 
$$
\Big\vert \fc{\mathcal J(S_F, b_0)}{M(J)}\Big\vert^{m+1}   
\le \max_{b_0} \big\vert \beta_d(b_0) \big\vert \ 
\left(\fc{m + 1}{2} |b_m|\right)^d \le 1,
$$ 
or
 \be\label{4.17}
\Big\vert \fc{\mathcal J(S_F, b_0)}{M(J)} \Big\vert \le \max_{b_0}  
\big\vert \beta_d(b_0)\big\vert^{1/(m+1)}  
\ \left(\fc{m + 1}{2} |b_m|\right)^{d/(m+1)} \le 1.  
\end{equation}  
In the case of an extremal function $F_{f_0}(z) = z + b_0^0 + b_m^0
z^{-m} + \dots$ for $\wt J(F)$, it must be 
$|\mathcal J(S_{F_{f_0}})/M(J)| = 1$, and (4.17) implies
$$
\fc{m + 1}{2} |b_m^0| = 1.
$$

Since the functions (4.16) with $m > 1$ satisfy $b_1 = \dots =
b_{m-1} = 0$, one can apply the well-known coefficient estimates of
Golusin and Jenkins (see [Go, Ch. XI], \cite{Je}) which provide in
our case the bound
  \be\label{4.18}
|b_m| \le 2/(m + 1)
\end{equation}
with equality only for $F = F_{m,t}$ with $|t| = 1$ (up to
translation $F_{m,t}(z) + c$). In this case, $M(J) =
|J(\kp_{m,\theta})|$.

We have established that any extremal function $f_0$ maximizing
$|J(f)|$ must be of the form (2.1). The theorem is proved.

\bk

\section{Proof of Theorem 2.2}

Note that from (1.3),
$$
a_n^2 - a_{2n-1} = b_1 b_0^{2n-4} + \text{lower terms with respect
to} \ b_0.
$$
We have to show that for all $m > 1$,
 \be\label{5.1}
|J_n(\kp_{m,\theta})| < J_n(\kp_0), \quad \kp_0 = z/(1 - z)^2.
\end{equation}
Then Theorem 2.1 implies that only the Koebe function is extremal
for Zalcman's functional.

This inequality is trivial for $m = 2$, because the series (1.7)
yields
$$
|J_n(\kp_{2,\theta})| \le 2 < J_n(\kp_0).
$$
For $m \ge 3$, we apply a result of \cite{Kr5} solving the
coefficient problem for univalent functions with \qc \ extensions
having small dilatations. Denote by $S_\iy(k)$ the subclass of $S^0$
consisting of the functions $f \in S$ having $k^\prime$-\qc \
extensions $\wh f$ to $\hC \ (k^\prime \le k)$ which satisfy $\wh
f(\iy) = \iy$, and let
$$
f_{1,t}(z) = \fc{z}{(1 - kt z)^2}, \quad |z| < 1, \ \ |t| = 1.
$$

\begin{prop} \cite{Kr5} For all
$f(z) = z + \sum\limits_2^\iy a_n z^n \in S_\iy(k)$ and all $k \le
1/(n^2 + 1)$, we have the sharp bound
 \be\label{5.2}
|a_n| \le 2 k/(n - 1),
\end{equation}
with equality only for the functions
 \be\label{5.3}
f_{n-1,t}(z) = f_{1,t}(z^{n-1})^{1/(n-1)} = z + \fc{2 k t}{n - 1}
z^n + \dots, \quad n = 3, 4, \dots
\end{equation}
\end{prop}

Note that every function (5.3) admits a \qc \ extension $\wh
f_{n-1,t}$ onto $\D^*$ with \Be \ coefficient $\mu_n(z) = t
|z|^{n+1}/z^{n+1}$ and $\wh f_{n-1,t}(\iy) = \iy$. Accordingly, $\wh
F_{n-1,t}(z) = 1/\wh f_{n-1,t}(1/z) \in \Sigma^0$ admits a \qc \
extension onto the unit disk with $\wh F_{n-1,t}(0) = 0$ and
$\mu_{\wh F_{n-1,t}}(z) = t |z|^{n-1}/z^{n-1}$ for $|z| < 1$.
Another essential point is that for any function
$$
F_{n-1}(z) := \wh F_{n-1,1}(z) = 1/\kp_0(1/z^{n-1})^{1/(n-1)},
$$
its homotopy disk $\{S_{F_{n-1}}\}$ in $\T$ is \Te \ geodesic.
Together with estimate (5.2), this implies that for any $m
> 2$ and small $r > 0$,
$$
|J_n(\kp_{m,r})| < r (n - 1 )^2;
$$
thus
$$
|J_n(\kp_{m,\theta})| < (n - 1 )^2 = J_n(\kp_0),
$$
completing the proof of Theorem 2.2.

\bk \noindent \textbf{Remarks}.

\medskip\noindent{\bf 1}. The above arguments work well also in the
case of functionals obtained by suitable perturbation of $J_n(f)$.
For example, one can take
 \be\label{5.4}
J(f) = a_n^2 - a_{2n-1} + P(a_3, \dots, \ a_{2n-2}),
\end{equation}
where $P$ is a homogeneous polynomial of degree $2n-2$,
$$
P(a_3, \dots, \ a_{2n-2}) = \sum\limits_{|k|=2n-2} c_{k_3, \dots,
k_n} a_3^{k_2} \dots a_n^{k_{2n-2}},
$$
and $|k| := k_3 + \dots \ + k_{2n-2}, \ a_j = a_j(f)$, assuming that
this polynomial has nonnegative coefficients and satisfies
$$
\max_S |P(a_3, \dots, \ a_{2n-2})| < \fc{(n - 1)^2}{2}.
$$
For any such functional, only the Koebe function is extremal.

\bk\noindent{\bf 2}. One can simplify the above proof applying
instead of Proposition 5.1 a weaker estimate $|a_n| \le 2 k/(n -
1)^2 + O(k^2)$ following from the variational formula for 
$f \in S(k)$.

\bk\noindent{\bf 3}. The evaluation of the coefficient functionals 
$J(F) = J(b_{m_1}(F), \dots \ , b_{m_p}(F))$   
on the class $\Sigma$ is somewhat different. 
In view of normalization, this case relates to the space the 
\uTs \ $\T$ (instead of $\T_1$ for the class $S$).  

The arguments exploited in the last step of the proof of Theorem 2.1 
do not work for the generic functionals on $\Sigma$, for example, 
if  $J(F) = b_m + \xi b_1$ with small $\xi$, because there are functions 
$F \in \Sigma$ for which the inequality (4.18) does not hold.

\bk

\section{Some new distortion theorems for higher coefficients}

As was mentioned, Theorem 3.1 provides various new distortion
estimates. For example, we obtain the following generalizations of
the inequality $|a_2^2 - a_3| \le 1$ to higher coefficients.

\begin{thm} For all $f \in S$ and integers $n > 3$
and $p \ge 1$,
$$
|a_n^p - a_2^{p(n-1)}| \le 2^{p(n-1)} - n^p.
$$
This bound is sharp, and the equality only occurs for the Koebe
function $\kp_\theta$.
\end{thm}

\medskip
\noindent \textbf{Proof}. Since $b_0 = - a_2$, the relation (3.1)
yields
$$
I_n(f) := a_n - a_2^{n-1} = (n - 2)(- 1)^{n-1} b_1 b_0^{n-3} +
\text{lower terms with respect to} \ \ b_0.
$$
This functional and $I_n^p(f) = |a_n^p - a_2^{p(n-1)}$ satisfy 
the assumptions of Theorem 2.1. The same
arguments as in the proof of Theorem 2.2 imply
$$
|I_n^p(\kp_{m,\theta})| < |I_n^p(\kp_\theta)| \quad \text{for all} \ \ m
\ge 2,
$$
completing the proof.

In the same way, one obtains

\begin{thm} For all $f \in S$ and integers $n > 2$ and $p \ge 1$,
$$
|a_{n+1}^p - a_2^p a_n^p| \le 2^p n^p - (n + 1)^p,
$$
with equality only for $f = \kp_\theta$.
\end{thm}

\section{Asymptotic theorems}

Another consequence of Theorem 2.1 concerns the asymptotic rate of
growth of generic homogeneous functionals on individual functions
with \qc \ extension (which are dense in $S$). We present here
somewhat restricted results.

\begin{thm} Let $\{J_j(f)\}_1^\iy$ be a sequence of uniformly bounded 
homogeneous
\hol \ (not necessarily distinct) functionals on $S$ of degrees $d_j
= d(J_j)$ satisfying the assumptions of Theorem 2.1 and such that
  \be\label{7.1}
|J_j(\kp_{m,\theta})| \le |J_j(\kp_0)| \quad \text{for all} \ \ m > 1. 
 \end{equation}
Then for any $f \in S^0$,
 \be\label{7.2}
 \limsup\limits_{j\to \iy} \ \Big\vert\fc{J_j(f)}{J_j(\kp_0)}\Big\vert
 = v(f) < 1.
 \end{equation}
\end{thm}

A similar result holds for the upper envelop $\sup_\a |J_\a (f)|$ of 
a family of uniformly bounded homogeneous \hol \ functionals
on $S$. 

\bk{\bf Proof}. It follows from (7.1) and Theorem 2.1 that only the
Koebe function $\kp_\theta$ is extremal for each $J_j$. We lift
$J_j(f)$ to \hol \ functionals $\wh J_j(\vp)$ on the space $\T_1$ 
(taking again $\vp = (S_f, - a_2(f))$) and consider the ratios  
$$
v(\vp) = \limsup\limits_{j\to \iy} \ \Big\vert 
\fc{\wh J_j(\vp)}{\wh J_j(\vp_0)} \Big\vert 
$$ 
where $\vp_0 = (S_{\kp_0}, - 2)$.  
The function $v(\vp)$ is well defined on this space and $v(\vp) \le 1$. 
Its upper semicontinued regularization 
$v^*(\vp) = \limsup\limits_{\vp^\prime \to \vp} v(\vp^\prime) \le 1$ 
is \psh \ on $\T_1$, hence by the maximum
principle it cannot attain the value $1$ inside $\T_1$;   
otherwise this function must be identically equal to $1$. But,
for example, $v^*(\mathbf 0) = v(\mathbf 0) = 0$, since near the
origin by Schwarz's lemma,
$$
\Big\vert \fc{\wh J_j(\vp)}{\wh J_j(\vp_0)}\Big\vert \le \const
\|\vp\|
$$
for all $p$. Therefore, $v(\vp) \le v^*(\vp) < 1$, which completes
the proof of (7.2).

\bk In the case of Zalcman's functional this yields that for any 
$f \in S^0$,
$$
v(f) = \limsup\limits_{n\to \iy} \fc{|a_n^2 - a_{2n-1}|}{(n - 1)^2}
< 1 
$$ 
(cf. \cite{Ha}, \cite{EV}). 
On the other hand, the bound (1.5) implies that $v(f) \le 1$ for
any $f \in S$. Similar estimates hold for perturbations of this
functional via (5.4).

\medskip
{\small\em{ \leftline{Department of Mathematics, Bar-Ilan
University} \leftline{52900 Ramat-Gan, Israel} \leftline{and
Department of Mathematics, University of Virginia,}
\leftline{Charlottesville, VA 22904-4137, USA}}}

\end{document}